
\input amstex
\magnification=1200


\def\skipfirstword#1 {}

\def\ir#1{\csname #1\endcsname}

\def\irrnSection#1#2{\edef\tttempcs{\ir{#2}}\head 
\expandafter\skipfirstword\tttempcs. #1 \endhead}

\def\irSubsection#1#2{\edef\tttempcs{\ir{#2}}\subhead 
\expandafter\skipfirstword\tttempcs. #1 \endsubhead}

\def\irSubsubsection#1#2{\edef\tttempcs{\ir{#2}}\subsubhead 
\expandafter\skipfirstword\tttempcs. #1 \endsubsubhead}

\long\def\prclm#1#2#3{\proclaim{\ir{#2}} #3 \endproclaim}

\def\Prf{\demo{Proof}}

\def\rem#1#2{\remark{{\bf \ir{#2} }} }

\def\bibitem#1{\ref\key #1 \by }

\long\def\eatit#1{}


\documentstyle{amsppt}

\pagewidth{5.5 in} 
\pageheight{7.5 in} 

\edef\tempqed{\qed}
\def\qed{\tempqed\enddemo}

\def\cal#1{\Cal #1}

\def\trans{H1}
\def\annalen{H2}
\def\duke{H3}
\def\bowdoin{H4}
\def\sundancetwo{H5}
\def\mtnwest{H6}
\def\Artin{Lemma II.9}

\def\C#1{\hbox{$\cal #1$}}
\def\pr#1{\hbox{{\bf P}${}^{#1}$}}

\topmatter 

\title Anticanonical Rational Surfaces\endtitle

\author Brian Harbourne\endauthor

\address
Department of Mathematics and Statistics
University of Nebraska-Lincoln
Lincoln, NE 68588-0323
\endaddress

\email bharbourne\@unl.edu, http://www.math.unl.edu/$\sim$bharbour/\endemail

\date January 24, 1996\enddate

\thanks This work was supported both by the National Science Foundation and
by a Spring 1994 University of Nebraska Faculty Development Leave.
I would also like to thank Tony Geramita for a helpful discussion,
and the referee for a careful reading of the paper. \endthanks

\subjclass Primary 14C20, 14J26; Secondary 14M20, 14N05\endsubjclass

\keywords Anticanonical, rational, surface, 
base points, fixed components, linear systems\endkeywords

\abstract A determination of the fixed 
components, base points and irregularity
is made for arbitrary numerically effective divisors on any smooth projective 
rational surface 
having an effective anticanonical divisor. 
All of the results are proven over an algebraically closed field
of arbitrary characteristic. Applications, to be treated in separate papers, 
include questions involving: points in good position,
birational models of rational surfaces in projective space, 
and resolutions for 0-dimensional subschemes
of $\pr2$ defined by complete ideals.\endabstract

\rightheadtext{Anticanonical Rational Surfaces}

\endtopmatter

\irrnSection{Introduction}{aaa}
This paper is concerned with complete
linear systems on smooth projective rational 
surfaces, over an algebraically closed field of arbitrary characteristic.
The focus is on {\it anticanonical\/} rational surfaces, i.e.,
those rational surfaces supporting an effective anticanonical divisor.
Such surfaces include all Del Pezzo surfaces, 
all blowings up of relatively minimal models of rational surfaces at
8 or fewer points, and all smooth complete toric surfaces,
but also include surfaces for which there is an effective 
but highly nonreduced anticanonical divisor; indeed, there is no bound 
to the multiplicities of fixed components of effective anticanonical divisors 
on anticanonical surfaces.

Anticanonical rational surfaces have 
received attention from numerous authors [Co], [D], [F], [\trans],
[\annalen], [\sundancetwo], [M], [Sk], [U] for various reasons, and share much
of the behavior of $K3$ surfaces [\bowdoin], [Ma], [PS], [SD], [St]. However, 
partly because the rational surfaces need not be relatively
minimal, and partly because their anticanonical classes are nontrivial,
the behavior of the rational surfaces is more complicated.
Thus previous work has been carried out under one or another special hypothesis,
to tame this anticanonical refractoriness. For example, [D] considers
blowings up of \pr 2 at $9$ or fewer points which are {\it en position presque g\'en\'eral},
and [F] assumes that there is a reduced anticanonical divisor,
and only considers classes of 
effective and numerically effective divisors which have no 
fixed component in common with the anticanonical divisor;
[\trans] handles arbitrary classes, but
considers only blowings up of \pr 2 for which there is a reduced and irreducible
anticanonical divisor. Likewise, [Sk]
mainly considers linear systems related to the anticanonical divisor itself.
In this paper we make no assumptions beyond the existence of an effective anticanonical divisor.

The structure of this paper is as follows.
Section II gathers results about rational surfaces generally, without
necessarily assuming effectivity of an anticanonical divisor. In Section III
we prove our main result,
Theorem III.1, which fully determines the behavior (that is,
dimensions, base points and fixed components) of complete linear systems 
for numerically effective classes
on smooth projective rational surfaces, under no special hypotheses beyond
effectivity of an anticanonical divisor. Our approach is
partly that of [\trans], and partly a modification of
[F] using an idea based on [A] (see \Artin). 

As a consequence of results in Section III,
we have for example the following theorem.
(Given a divisor class \C L
on a surface $X$, we recall that
a class \C L is {\it numerically effective\/} if it meets every effective
divisor nonnegatively. We will also find it convenient to
identify an invertible sheaf with its corresponding divisor class, and 
employ the additive notation customary for the group of divisor classes.
We will usually reserve $\otimes$ to denote restriction, as in the restriction
$\C O_C\otimes\C F$ of a divisor class \C F on a surface $X$ to a curve $C\subset X$.)

\prclm{Theorem}{intr}{Let \C F be a numerically effective divisor class on
a smooth projective rational anticanonical surface $X$ and let $D$ be a general section
of $-K_X$. Then $h^0(X, \C F)>0$; moreover, $h^1(X, \C F)>0$ if and only if
$\C F\cdot K_X=0$ and a general section of $\C F$ has a connected 
component disjoint from $D$. In fact, $1+h^1(X, \C F)$ is the number of connected
components of a general section of $\C F-K_X$.}

A proof is given at the end of Section III, but a remark here may be helpful.
Although $h^0(X, \C F)>0$ is considered in Section II,
the main interest is in $h^1$. That there be a topological sufficient
condition for $h^1$ to be positive is easy to see.
If a section of $\C F$ has a connected 
component disjoint from $D$, 
then $\C F-K_X$ has a section $C$ which is not
connected, so $h^0(C, \C O_C)>1$. 
But $h^1(X, -(\C F-K_X))=h^1(X, \C F)$ by duality, so
from $0\to -(\C F-K_X)\to \C O_X\to \C O_C\to 0$,
we see $h^1(X, \C F)=-1+h^0(C, \C O_C)$ is positive.

It is more difficult to see that the question is
purely topological, since $C$ need not be reduced,
hence {\it a priori}, 
$h^0(C, \C O_C)$ could exceed the number of connected components
of $C$. That it does not is a consequence
of the detailed analysis of the fixed components of
numerically effective divisors on anticanonical rational surfaces,
carried out in Section III, showing, among other things, 
that the fixed part $N$ of the linear system of sections of
a numerically effective class is reduced and no
components of $N$ can be components of an anticanonical divisor, unless $N$ contains
an entire anticanonical divisor. 

Some readers may be interested in classes that are not necessarily numerically
effective; a word to them may be in order. Given an arbitrary class \C F
on a smooth projective rational anticanonical surface $X$,
to determine $h^0(X, \C F)$, and fixed components and base points
of the linear system of sections of \C F
when \C F is the class of an effective divisor, one essentially needs to know
the monoid EFF of classes of effective divisors. However, this is reasonably
accessible and is discussed in more detail in the remark at the end of
Section III.

Applications, to questions involving points in good position,
to birational models of rational surfaces in projective space, 
and to computing resolutions for 0-dimensional subschemes
of \pr2 defined by complete ideals will be
considered in separate papers.

\irrnSection{Results not assuming $-K$ effective}{bbb}
In this section we consider smooth 
projective rational surfaces, stating some results which hold in general,
whether the surfaces are anticanonical or not. We begin with the homomorphism 
on Picard groups induced by a morphism of schemes. In the cases of 
interest here, this homomorphism
is well-known to behave very nicely. (Given 
a surface $X$, a curve $C\subset X$, and a divisor class \C F on $X$,
we may in place of the more accurate
but more complicated $h^i(C, \C F\otimes\C O_C)$ write
$h^i(C, \C F)$ for the dimension of the $i^{\hbox{th}}$ cohomology group
of the restriction of \C F to $C$.)

\prclm{Lemma}{one}{Let $\pi:Y\to X$ be a birational morphism of 
smooth projective surfaces,
$\pi^*:\hbox{Pic}(X)\to\hbox{Pic}(Y)$ the corresponding homomorphism 
on Picard groups,
and let \C L be a divisor class on $X$.
\itemitem{(a)} The map $\pi^*$ is an injective intersection-form preserving
map (of free abelian groups of finite rank if $X$ is rational).
\itemitem{(b)} The map $\pi^*$ preserves dimensions of cohomology groups;
i.e., $h^i(X,\C L)=h^i(Y,\pi^*\C L)$ for every $i$.
\itemitem{(c)} The map $\pi^*$ preserves EFF;
i.e., \C L is the class of an effective divisor if and only if $\pi^*\C L$ is.
\itemitem{(d)} The map $\pi^*$ preserves numerical effectivity;
i.e., $\C L\cdot F\ge 0$ for every effective divisor $F$ on $X$ if and only if 
$(\pi^*\C L)\cdot F'\ge 0$ for every effective divisor $F'$ on $Y$.}

\Prf (a) See [Ha, V]. 

(b) Use [Ha, V.3.4] and the Leray spectral sequence.

(c) This follows from (b) with $i=0$.

(d) If $\pi^*\C L$ is numerically effective, then numerical 
effectivity for 
\C L follows from (a) and (c). The converse follows from
[Ha, V.5] and induction.  \qed

Given an appropriate scheme $X$, we recall that $K_X$ 
denotes its canonical class.

\prclm{Lemma}{RR}{Let $X$ be a smooth projective 
rational surface, and let \C F be
a divisor class on $X$.
\itemitem{(a)} We have: $h^0(X, \C F)-h^1(X, \C F) + h^2(X, 
\C F)= (\C F^2-K_X\cdot \C F)/2+1$.
\itemitem{(b)} If \C F is the class of an effective divisor, then $h^2(X, \C F)=0$.
\itemitem{(c)} If \C F is numerically effective,
then $h^2(X, \C F)=0$ and $\C F^2\ge 0$.}

\Prf Item (a) is just the Riemann-Roch formula in the case of a 
rational surface. Item (b) follows by duality, while item (c)
is elementary. \qed

\prclm{Corollary}{NefIsEffCor}{On a smooth projective rational surface, 
a numerically effective divisor meeting the anticanonical 
class nonnegatively is in EFF. In particular,
effectivity of an anticanonical divisor implies effectivity of all numerically 
effective divisors.}

\Prf This follows from \ir{RR}(a) and (c). \qed

The next lemma is a well-known consequence of the Hodge index Theorem.

\prclm{Lemma}{Lattice}{Let $X$ be a smooth projective 
rational surface, and denote by
$K^\perp$ the subspace of Pic$(X)$ perpendicular to the canonical class $K_X$.
Then $K^\perp$ is negative definite if and only if $K_X^2>0$.
Also, $K^\perp$ is negative semidefinite if and only if $K_X^2=0$; in
this case, if $x\in K^\perp$, then $x^2=0$ if and only if $x$ is a multiple of
$K_X$. }

In spite of the title of this section, the next fact
involves anticanonical surfaces, but it is convenient 
to state it now.

\prclm{Lemma}{ODisirr}{Let $X$ be 
a smooth projective anticanonical rational surface. 
Let $D$ be an effective anticanonical divisor on $X$
and let \C F be any divisor class on $X$. Then
$h^0(D, \C O_D)=h^1(D, O_D)=1$, and 
$h^0(D, \C F)-h^1(D, \C F)=-K_X\cdot \C F$.}

\Prf By duality, $h^2(X, K_X)=h^0(X, \C O_X)=1$ and
$h^1(X, K_X)=h^1(X, \C O_X)$ vanish, while
$h^0(X, K_X)=0$ since $X$ is rational. 
Thus from $0\to K_X\to \C O_X\to \C O_D\to 0$ we see that
$h^0(D, \C O_D)=1=h^1(D, \C O_D)$.
From $0\to \C F+K_X\to \C F\to \C F\otimes\C O_D\to 0$, comparing
$h^0-h^1$ for $\C F\otimes\C O_D$ against $h^0-h^1+h^2$ for 
\C F and $\C F+K_X$, using Riemann-Roch for the latter two
and simplifying,  we see that 
$h^0(D, \C F)-h^1(D, \C F)=-K_X\cdot \C F$. \qed

We will need a Bertini-type theorem. Recall that a fixed 
component free linear system is
said to be {\it composed with a pencil}, if away from the base points
it defines a morphism whose image has dimension $1$.

\prclm{Lemma}{Bertini}{Let $X$ be a smooth projective rational surface,
and let \C F be the class of a nontrivial effective divisor on $X$ without fixed components.
\itemitem{(a)}If the linear system of sections of \C F is composed with a pencil, then there is an $r>0$
such that every section of \C F is the divisorial
sum of $r$ sections of some class \C C whose 
general section is reduced and irreducible 
and whose sections comprise a pencil. 
Moreover, if $\C F\cdot K_X\le 0$, then either: $r=1$, $\C F^2=2$,
$\C F\cdot K_X=0$, $h^1(X, \C F)=0$ and $-K_X$ is not effective; 
$r=\C F^2=-K_X\cdot\C F=1$ and $h^1(X, \C F)=0$;  
$\C F^2=0$, $-K_X\cdot\C F=2r$ and $h^1(X, \C F)=0$; or
$\C F^2=K_X\cdot\C F=0$ and $h^1(X, \C F)=r$.
\itemitem{(b)}If the linear system of sections of \C F is not composed with a pencil, 
then $\C F^2>0$ and the general 
section of \C F is reduced and irreducible.}

\Prf We first point out that the general section of 
the class \C F of a fixed component free divisor is
always reduced. For if not, then by
[J, Corollaire 6.4.2] the morphism defined by the 
sections away from the base points factors through
Frobenius. This cannot happen for complete linear systems
on rational surfaces, since it would imply that
the characteristic $p$ is nonzero, that $\C F=p\C G$ for some class
\C G, and (since $X$ is rational) that $h^0(X, \C G)=h^0(X, p\C G)$. 
Since the sections of \C F are fixed component free,  
$h^0(X, \C G)=h^0(X, p\C G)>1$,
hence $h^0(X, p\C G)\ge \hbox{dim Sym}^p(H^0(X, \C G))=
\left(\matrix d-1+p\\
                p\endmatrix\right)>d$, where $d=h^0(X, \C G)$.

For (b), note that $\C F^2=0$ would imply the linear system of sections of 
\C F is composed with a pencil; for the irreducibility,
use [J, Th\'eor\`eme 6.3].

Now consider (a). Since $X$ is rational, the linear system of sections of \C F is composed with
a rational pencil.
By appropriately blowing up the base points of the sections of \C F,
we obtain a class $\C F'$ on a surface
$Z$ whose sections are base point free, whose general section is the proper transform
of a general section of \C F. 
The sections of $\C F'$ define a morphism $Z\to R$, where $R$ is isomorphic
to \pr 1 since the pencil is rational. By Stein factorization [Ha, III.11.5],
there is a morphism $Z\to \pr 1$, with connected fibers, factoring $Z\to R$.
If $r$ is the degree of $\pr 1\to R$ and $\C C'$ is the class of
a fiber of $Z\to \pr 1$, then $\C F'=r\C C'$, where the sections of $\C C'$
are connected. Taking \C C to be the class
of a curve on $X$ whose proper transform is a general section of
$\C C'$, we see that $\C F=r\C C$, and hence (comparing sections of
\C F with sections of $\C F'$)
that every section of \C F is the divisorial
sum of $r$ sections of \C C. Thus it is enough to check that a general
section of $\C C'$ is integral, but this is well-known to hold
(cf. [J, Th\'eor\`eme 4.10] and [J, Th\'eor\`eme 3.11])
since $Z\to \pr 1$ induces a seperable algebraically closed extension
of function fields.
This proves the first part of (a).

For the rest, assume $\C F\cdot K_X\le 0$. 
Because \C C gives a pencil, $h^0(X, \C C)=2$. Because the sections of \C F are 
composed with this pencil, $h^0(X, \C F)=h^0(X, r\C C)=r+1$. By Riemann-Roch,
$r+1=h^0(X, \C F)\ge (r^2\C C^2-r\C C\cdot K_X)/2+1$. Now suppose
\C C has positive self-intersection. If $\C C\cdot K_X=0$, then $\C C^2$
must be even and thus is at least $2$,
so $r+1\ge r^2\C C^2/2+1\ge r^2+1$, hence $r=1$ and $\C C^2=2$. 
Also, from Riemann-Roch we can now calculate that $h^1(X, \C F)=0$.
But if $-K_X$ were effective, then, since the sections of \C C are fixed component free, 
$-K_X\cdot \C C=0$ means $D$ is disjoint from a general section of \C C,
where $D$ is some effective anticanonical divisor. 
Thus $\C O_D\otimes\C C=\C O_D$.
By \ir{ODisirr}, $h^1(D, \C O_D)=1$. Since  
$h^2(X, K_X+\C C)=h^0(X, -\C C)=0$, from 
$0\to K_X+\C C\to \C C\to \C O_D\otimes \C C\to 0$ 
we get the contradiction that $h^1(X, \C C)>0$. 
Thus $-K_X$ cannot be effective.

If, on the other hand, $\C C^2>0$ and $\C C\cdot K_X<0$, then 
$2r\ge r^2+r$ by Riemann-Roch so $r=1$ (and $\C F=\C C$) whence 
$2=(\C C^2-\C C\cdot K_X)/2 +1+h^1(X, \C C)=2+h^1(X, \C C)$, and 
thus $\C C^2=\C C\cdot (-K_X)=1$ and $h^1(X, \C C)=0$.

Finally, say $\C C^2=0$. Then $-K_X\cdot \C C$ is even, but
$2r\ge -r\C C\cdot K_X$ by Riemann-Roch, so $-K_X\cdot \C C$
is either $2$ or $0$. From Riemann-Roch we have in the first case
$h^0(X, \C F)=r+1+h^1(X, \C F)$, and we have in the second 
$h^0(X, \C F)=1+h^1(X, \C F)$, so from $r+1=h^0(X, \C F)$
follow $h^1(X, \C F)=0$ and $h^1(X, \C F)=r$, respectively.  \qed

\prclm{Lemma}{fcfreeandreg}{Let $X$ be a smooth projective rational surface,
and let \C C be the class of an effective divisor without fixed components. If
$\C C\cdot K_X<0$, then $h^1(X, \C C)=0$.}

\Prf Say \C C is the class of a reduced and irreducible divisor $C$ on $X$,
with $C\cdot K_X<0$. Then $h^1(C, \C O_X(rC)\otimes\C O_C)=0$ for all $r>0$.
(If $C$ is smooth this follows by duality, since by adjunction
$\C O_X(-rC)\otimes K_C$ has negative degree. In general,
$h^1(C,\C O_C)=h^2(X, \C O_X(-C))=1+(C^2+C\cdot K_X)/2$ 
follows from checking the cohomology groups for 
the sequence $0\to \C O_X(-C)\to \C O_X\to \C O_C\to 0$, and, for the
second equality, applying Riemann-Roch. Thus $\C O_X(rC)\otimes\C O_C$
has degree greater than $2h^1(C,\C O_C)-2$;  
$h^1(C, \C O_X(rC)\otimes\C O_C)=0$ for all $r>0$ now follows by Proposition 7
on page 59 of [D].) Since $X$ is rational, $h^1(X, \C O_X)=0$. Induction 
on $r$ applied to the cohomology groups of 
$0\to \C O_X((r-1)C)\to \C O_X(rC)\to \C O_C(rC)\to 0$ now gives the result.

In general, if the linear system of sections of \C C
has no reduced and irreducible member,
then, by \ir{Bertini}, we find that \C C is the class of a multiple
of a reduced irreducible divisor $C'$ which 
moves in a pencil at least. The result now 
follows as above for any positive multiple
of $C'$ and hence for any positive multiple of $C$. \qed

The following lemma, essentially Theorem 1.7 of [A], 
turns out to be the key to avoiding special assumptions (such as are
assumed in [F])
regarding components in common with an anticanonical divisor
when handling numerically effective classes
on anticanonical surfaces.   

\prclm{Lemma}{keylemma}{Let $X$ be a smooth projective surface supporting
an effective divisor $N$ and a divisor class $\C F$ which meets every 
component of $N$ nonnegatively. If $h^1(N, \C O_N)=0$, 
then $h^0(N, \C F)>0$ and $h^1(N, \C F)=0$.}

\Prf Modifying the proof of [A, Theorem 1.7]
gives the result. See Lemma II.6 of [\mtnwest] for an explicit proof. \qed

\prclm{Corollary}{keycor}{Let $X$ be a smooth projective 
rational surface 
and let \C N be the class of a nontrivial effective divisor $N$ 
on $X$. If $\C N+K_X$ is
not the class of an effective divisor and \C F meets every component of $N$ nonnegatively
(in particular, if \C F is numerically effective), 
then $h^0(N, \C F)>0$, $h^1(N, \C F)=0$, $N^2+N\cdot K_X<-1$,
and every component $M$ of $N$ is a smooth rational curve (of negative
self-intersection, if $M$ does not move).}

\Prf Clearly, $h^0(X, -\C N)=0$, while 
$h^2(X,-\C N)=h^0(X,\C N+K_X)=0$ 
follows by duality and by hypothesis. By Riemann-Roch,
$-2\ge-2-2h^1(X, -\C N)=N^2+N\cdot K_X$. 
Next, taking cohomology of $0\to -\C N\to \C O_X\to \C O_N\to 0$,
we see $h^1(N, \C O_N)=0$. Since 
\C F meets every component of $N$ nonnegatively, applying \ir{keylemma}
shows $h^0(N, \C F)>0$ and $h^1(N, \C F)=0$. To finish,
apply the foregoing to a component $M$ of $N$.
Then $M^2+M\cdot K_X=-2$, so $M$ is smooth 
and rational. Substituting into Riemann-Roch 
gives $h^0(X, \C O_X(M))\ge M^2+2$, hence, if $M$ does not move, then
$1\ge M^2 +2$ so $M^2\le -1$. \qed

\irrnSection{Results assuming $-K$ effective}{ccc}
We now develop results concerning anticanonical rational surfaces.
The results of this section give a complete determination
of $h^1(X, \C F)$, the class \C N of the fixed part of 
the linear system of sections of \C F
and the base points of the linear system of sections of 
\C F (note that $\C F\in\hbox{EFF}$ is automatic 
by \ir{NefIsEffCor}), for any numerically effective class
\C F on a smooth projective rational anticanonical surface $X$.

Briefly, what we find is this. Behavior of the linear system
of sections of a numerically effective class \C F is controlled by
$-K_X\cdot \C F$. If this is at least 2, then \C F is regular
and its linear system has empty base locus. 

If $-K_X\cdot \C F=1$,
\C F is still regular, but its linear system has a base point,
and possibly a base divisor. (The base divisor is always the class of
a reduced but not necessarily irreducible exceptional curve.
The pencil of cubics through eight general points gives an example
of a numerically effective class \C F on a blow up of \pr2 at eight points
whose sections have a base point; by blowing up the base point and
pulling back to the blown up surface, the base point converts to
a divisorial base locus consisting of the exceptional curve $E$
of the last blow up. In this example, $\C F\cdot E=0$; this example is not
completely general, since
it is also possible to have a divisorial base locus $E$ with
$\C F\cdot E>0$, as happens if $X$ is the blow up of
\pr2 at the base points of a general pencil of cubics
and \C F is $-rK_X+\C E$, where \C E is the class of the exceptional curve
$E$ coming from one of the blow ups.)

Regardless of what else happens, if $-K_X\cdot \C F=1$, then the sections
of \C F have a base point at a smooth point of any 
effective anticanonical divisor $D$. Blowing this base point up and considering
$\C G=\C F-\C E'$ on the blown up surface (where $\C E'$ is the class of the exceptional divisor
of the blow up) shows that numerically effective classes \C G can occur
with $-K_X\cdot \C G=0$, such that $h^1(X,\C G)>0$ and either the sections of
\C G are base point free (which happens if the sections of \C F 
are fixed component free) or the base locus consists of a smooth
rational curve $N$ with $N^2=-2$ (which happens, for example,
if we take $\C F=-rK_X+\C E$ with $r>1$, as in the preceding paragraph;
note that $r=1$ does not work since then \C G is not numerically effective).
It turns out that the only other possibility that can occur for a numerically effective class
perpendicular to $K_X$ is that the base locus of its linear system
of sections can contain an anticanonical divisor, and understanding that this may happen
in essentially only one way (this being that of a numerically effective class
perpendicular to $-K_X$ whose restriction to an effective anticanonical divisor is
nontrivial---although a precise description of the base locus in this case
depends on whether or not the numerically effective class has 
self-intersection 0)
completes our understanding of linear systems of numerically effective classes
on anticanonical surfaces.

The following theorem gives a rigorous statement of the preceding discussion,
summarizing our results in this section:

\prclm{Theorem}{mainthm}{Let $X$ be a smooth 
projective rational anticanonical surface
with a numerically effective class \C F and 
let $D$ be a nonzero section of $-K_X$. 
Let \C N be the class of the
fixed part of the linear system of sections of \C F,
and let $\C H=\C F-\C N$ be the class of the free part.
\itemitem{(a)} If $-K_X\cdot\C F\ge 2$, then 
$h^1(X, \C F)=0$ and the sections of \C F are base point (and thus
fixed component) free (i.e., $\C N=0$).
\itemitem{(b)} If $-K_X\cdot\C F=1$, then $h^1(X, \C F)=0$. 
If the sections of \C F are fixed component
free, then the sections of \C F have a unique base point, which is on $D$.
Moreover, the linear system of sections of \C F has a fixed component if and only if
$\C H=r\C C$ and $\C N=\C N_1+\cdots+\C N_t$,
where $\C C\in K_X^\perp$ is a class with
$h^1(X, \C C)=1$ whose general section
is reduced and irreducible,
$r=h^1(X, \C H)$ with $r>1$ only if $\C C^2=0$, $\C N_i$ 
is a smooth rational curve for every $i$, $\C N_i^2=-2$ 
and $\C N_i\cdot \C N_{i+1}=1$ for $i<t$, $\C N_t^2=-1$, $\C N_i\cdot \C N_j=0$
for $j>i+1$, $\C C\cdot\C N_1=1$, 
and $\C C\cdot\C N_i=0$ for $i>1$.
\itemitem{(c)} If $-K_X\cdot\C F=0$, then either $\C N=0$ (in which case
the sections of \C F are base point free,
$\C F\otimes \C O_D$ is trivial and either 
$\C F^2>0$ and $h^1(X, \C F)=1$ or $\C F=r\C C$
and $h^1(X, \C F)=r$, where $r>0$ and \C C is 
a class of self-intersection $0$ whose general section
is reduced and irreducible), 
or \C N is a smooth rational curve
of self-intersection $-2$ (in which 
case $h^1(X, \C F)=1$, $\C N\otimes \C O_D$ is trivial,
and $\C H=r\C C$, where $r>1$ and \C C is reduced and irreducible with 
$\C C^2=0$, 
$\C C\cdot \C N=1$ and $\C C\otimes \C O_D$ being trivial), 
or $\C N+K_X\in\hbox{EFF}$.
\itemitem{(d)} We have $\C N+K_X\in\hbox{EFF}$ if and only if $\C F\cdot D=0$ but
$\C F\otimes \C O_D$ is nontrivial. In this case, there exists a birational
morphism of $X$ to a smooth
projective rational anticanonical surface $Y$, 
and either: $K_Y^2<0$, there is a
numerically effective class $\C F'$ on $Y$, \C F is the pullback of $\C F'-K_Y$, and
$0=h^1(Y, \C F')=h^1(X, \C F)$; 
or $K_Y^2=0$, \C H and \C N are the pullbacks
of $-sK_Y$ and $-rK_Y$ for some integers $s\ge 0$ and $r>0$ respectively, 
and $h^1(X, \C F)=\sigma$, where $\sigma=0$ if 
$s=0$, and otherwise $r<\tau$ and $\sigma=s/\tau$, where
$\tau$ is the least positive integer such that 
the restriction of $-\tau K_X$ to $D$ is trivial.}

At the end of this section we show how this summary statement
follows from the lemmas we will by then have accumulated. 

On a K3 surface, one never needs to worry about a linear system's having
a fixed component which is also a component of an effective
anticanonical divisor. As a consequence, the analysis of 
complete linear systems on K3 surfaces is less complicated than
on rational surfaces with an effective anticanonical divisor.
That this analysis can, nevertheless, be carried out
for rational surfaces is perhaps partly due to the occurrence of such
behavior being very restricted, as this next result shows:
 
\prclm{Corollary}{allornone}{Let $X$ be a smooth 
projective rational anticanonical surface
with a numerically effective class \C F and let \C N be the class of
the fixed part of the linear system of sections of \C F. Then either no
fixed component of the linear system of sections of \C F is a component of any section of
$-K_X$, or \C N contains an anticanonical divisor (i.e., $\C N+K_X\in\hbox{EFF}$).}

Thus this corollary implies that 
no numerically effective divisor on $X$ can have a fixed 
component in common with the anticanonical linear system if $h^0(X, -K_X)>1$
(and so in particular if $K_X^2>0$). 
Note that this corollary can be obtained from \ir{mainthm},
but for convenience we put off the proof until the end
of this section, obtaining it
as a consequence of other results later in this section.

The next corollary gives some geometrical insight into the occurrence
of irregular numerically effective classes; we give its proof now
as an example of applying \ir{mainthm}.

\prclm{Corollary}{regcond}{Let $X$ be a smooth 
projective rational anticanonical surface
with a numerically effective class \C F and let
$D$ be a nonzero section of $-K_X$. 
\itemitem{(a)} We have $h^1(X, \C F)>0$ 
if and only if the general section of \C F has a 
connected component disjoint from $D$. 
\itemitem{(b)} Say $\C F^2-K_X\cdot \C F>0$; then $h^1(X, \C F)>0$ 
if and only if the general section of \C F is disjoint from $D$.}

\Prf Let \C N denote the class of the fixed part of the linear system of sections of \C F,
and \C H the rest. Suppose
$h^1(X, \C F)>0$. By \ir{mainthm}, we see that $\C F\cdot D=0$.
If $\C N +K_X\notin\hbox{EFF}$, then the linear systems of sections of \C F and $-K_X$
have no common fixed components by \ir{allornone}, so $\C F\cdot D=0$
implies that the general section of \C F is disjoint from $D$.
If $\C N +K_X\in\hbox{EFF}$, then by \ir{mainthm}(d) $h^1(X, \C F)>0$
implies that \C H and \C N are the pullbacks of $-sK_Y$ and $-rK_Y$,
where $K_Y^2=0$, and thus that $\C F^2-K_X\cdot \C F=0$.
But \C H is numerically effective and its linear system of sections is
without fixed components, 
so $\C H\cdot \C N=-sK_Y\cdot(-rK_Y)=0$ and 
$\C H\cdot (-K_X)=-sK_Y\cdot (-K_Y)=0$
imply that \C H has a section containing a connected component of 
a general section of \C F disjoint from $D$. This establishes
the forward implications of (a) and (b). 

For the converses, note that if a general section 
of \C F has a connected component
disjoint from $D$, then some section $C$ of 
$\C F-K_X$ is not connected, hence $h^1(C, K_C)=h^0(C, \C O_C)>1$.
Now $h^1(X, \C F)>0$ follows from $K_C=\C F\otimes\C O_C$ and from
$0\to K_X\to \C F\to \C F\otimes \C O_C\to 0$.  \qed

We now begin to develop the tools we will need to prove \ir{mainthm}.
The following lemma comes from [F], included here for convenience,
and in order to clarify the statement of the lemma.

\prclm{Lemma}{RFa}{Let \C L be a divisor class on a smooth projective 
rational anticanonical surface $X$ with a reduced connected section $L$
no component of which is a fixed component of the linear system of sections of $-K_X$.
\itemitem{(a)} If $-K_X\cdot \C L>0$, then $h^1(X, \C L)=0$.
\itemitem{(b)} If $-K_X\cdot \C L=0$, then $h^1(X, \C L)=1$.}

\Prf (a) This is essentially [F, 5.2], although the 
statement there
does not make it clear that it is not enough
that $L$ be reduced with each of its connected 
components meeting
$-K_X$ positively. The following proof is taken from [F].

From $0\to \C O_X\to \C L\to \C O_L\otimes \C L\to 0$, 
we obtain $h^1(X, \C L)=h^1(L, \C L)$.
By adjunction, $\C O_L\otimes \C L=K_L\otimes(-K_X)$, so by duality
$h^1(L, \C L)=h^0(L, K_X)$. But the restriction of $-K_X$ to $L$ is 
the class of an effective divisor (since no component of $L$ is a fixed component of 
the linear system of sections of $-K_X$)
and nontrivial (since $-K_X$ meets $L$ positively); since $L$ is reduced
and connected, it follows that $h^0(L, K_X)=0$.

(b) This is [F, 5.3]. The proof is similar to (a), except 
here the hypotheses
imply that the restriction of $-K_X$ to $L$ is trivial, so now
$\C O_L\otimes \C L=K_L\otimes(-K_X)=K_L$. Thus $h^1(L, \C L)=h^0(L, 
\C O_L)$,
and $h^0(L, \C O_L)=1$ since $L$ is reduced and connected. \qed

The next lemma is a modification of Lemma 9 of [F], 
in which the requirement there 
that $-K_X$ have a section meeting $\C L_2$ in just a finite set of points 
is weakened here by only assuming 
that $h^0(X, \C L_2+K_X)=0$. Thus in the 
following lemma, the fixed part of the linear system of sections of
$\C L$ contains a section of $\C L_2$, and 
{\it a priori\/} could 
contain fixed components of 
the linear system of sections of 
$-K_X$, but merely cannot contain an entire anticanonical divisor.

\prclm{Lemma}{RFb}{Let $X$ be a smooth projective rational anticanonical 
surface.
Let $\C L_1$ and $\C L_2$ be classes of effective divisors on $X$ such that 
$h^0(X, \C L_2+K_X)=0$ and such that
$\C L=\C L_1+\C L_2$ and $\C L_1$ are numerically effective
with $h^0(X, \C L)=h^0(X, \C L_1)$ and 
$\C L_2\ne 0$; then $\C L_1\cdot \C L_2>0$.}

\Prf Suppose on the contrary that $\C L_1\cdot \C L_2=0$.
Since $\C L_1$ is numerically effective,
this means $\C L_1\cdot C=0$ for every component $C$ of 
the linear system of sections of $\C L_2$.  
Since $\C L=\C L_1+\C L_2$ is numerically effective,
it follows that $\C L_2\cdot C\ge 0$ for every component
$C$ of the linear system of sections of $\C L_2$, and hence that $\C L_2$ is
numerically effective; thus $\C L_2^2\ge 0$ and $-K_X\cdot \C L_2\ge 0$.
But $\C L_2$ is the class of a unique effective divisor, so from Riemann-Roch
we have $1\ge (\C L_2^2-K_X\cdot \C L_2)/2+1$, 
hence $\C L_2^2=-K_X\cdot \C L_2=0$, so $\C L_2=0$ by \ir{keycor}. \qed

\prclm{Lemma}{sizero}{Let \C F be a nontrivial numerically effective class on a
smooth projective rational anticanonical surface $X$, such that $\C F^2=0$ and
the class of the fixed part of the linear system of sections of \C F is \C N. 
If $h^0(X, \C N+K_X)=0$, then the sections of \C F are base point free
and composed with a pencil.}

\Prf Let $\C F-\C N=\C H$; then $0=\C F^2=\C H^2+\C H\cdot\C N
+\C N\cdot\C F$, hence by numerical effectivity of \C F and \C H, we have
$\C H\cdot \C N=0$, so $\C N=0$ by \ir{RFb}. Thus the sections of \C F are not only fixed
component free, but, since $\C F^2=0$, base point free and composed
with a pencil.  \qed

\prclm{Lemma}{siposirrzero}{Let \C F be a numerically effective class on a
smooth projective rational anticanonical surface $X$, such that $\C F^2>0$,
the class of the fixed part of the linear system of sections of \C F
is \C N and $\C F-\C N=\C H$. Let $D$ be a nonzero section of $-K_X$.
If $h^0(X, \C N+K_X)=0=h^1(X, \C H)$, then $\C F\cdot D>0$, \C F is regular, and
the linear system of sections of \C F is fixed component free.}

\Prf First we show $\C N=0$. Suppose not.
Let $N$ be the nontrivial section of \C N.
Then by \ir{keycor}, $h^1(N, \C F)=0$ and $h^0(N, \C F)>0$.
Using regularity of \C H we see
$0\to \C H\to \C F\to \C F\otimes\C O_N\to 0$
is exact on global sections, so $h^0(N, \C F)>0$
guarantees that \C H cannot be all of the class of the nonfixed part 
of \C F, contradiction. Thus \C N is trivial,
and the sections of $\C F=\C H$ are fixed component free, and regular by assumption. 
To see that $\C F\cdot D>0$, consider
$0\to \C F+K_X\to \C F\to \C F\otimes\C O_D\to 0$.
By duality, $h^2(X, \C F+K_X)=0$, so $h^1(X, \C F)\ge h^1(D, \C F\otimes\C O_D)$.
If $\C F\cdot D=0$, then the general section of \C F is disjoint from $D$,
so $\C F\otimes\C O_D=\C O_D$, but this contradicts regularity of \C F, since
$h^1(D, \C O_D)=1$ by \ir{ODisirr}. \qed

The argument for the next result is essentially that of [F, 10.5],
modified to fit the current more general situation.

\prclm{Lemma}{siposhirr}{Let \C F be a numerically effective class on a
smooth projective rational anticanonical surface $X$, with $\C F^2>0$,
where \C N is the class of the fixed part $N$ of the linear system of sections of \C F
and $\C F-\C N=\C H$ is the class of the free part. 
Let $D$ be an arbitrary anticanonical divisor on $X$.
Suppose $h^0(X, \C N+K_X)=0$ but $h^1(X, \C H)>0$. 
\itemitem{(a)} If $\C N=0$, then $\C F\cdot K_X=0$ and $h^1(X, \C F)=1$.
\itemitem{(b)} If $\C N\ne 0$, then $0\le -K_X\cdot\C F\le 1$, and there is a 
class $\C C\in K_X^\perp$ with $h^1(X, \C C)=1$ whose general section
is reduced and irreducible
such that $\C H=r\C C$, where $h^1(X, \C H)=r$.
Moreover, no component of $N$ is a component of $D$
and either: $h^1(X, \C F)=1$, in which case $\C F\cdot K_X=0$,
$\C C^2=0$, $r>1$, $\C N\cdot\C C=1$  and \C N is the class of a smooth
rational curve of self-intersection $-2$; 
or $h^1(X, \C F)=0$,
in which case $\C N=\C N_1+\cdots+\C N_t$,
where $\C N_i$ is a smooth rational curve for every $i$, $\C N_i^2=-2$ 
and $\C N_i\cdot \C N_{i+1}=1$ for $i<t$, $\C N_t^2=-1$,  $\C N_i\cdot \C N_j=0$
for $j>i+1$, $\C C\cdot\C N_1=1$, $\C C\cdot\C N_i=0$ for $i>1$,
and either $r=1$ or $\C C^2=0$.}

\Prf (a) By \ir{fcfreeandreg}, we see $\C F\cdot K_X=0$. 
By \ir{Bertini}, \C F is the class of a reduced and irreducible
divisor $F$, and the restriction of $K_X$ to $F$ is trivial,
since the sections of \C F are fixed component free, $-K_X$ is in EFF and 
$\C F\cdot K_X=0$.
From adjunction we thus have $K_F=\C F\otimes\C O_F$, so 
$h^1(F, \C F\otimes\C O_F)=1$. Now
$h^1(X, \C F)=1$ follows from 
$0\to \C O_X\to \C F\to \C F\otimes\C O_F\to 0$. 

(b) By \ir{fcfreeandreg}, we see $\C H\cdot K_X=0$,
and by \ir{RFb}, $\C H\cdot \C N>0$. 

Also, $\C H=r\C C$ and $h^1(X, \C H)=r$ for some 
positive integer $r$, where \C C is the class of some reduced 
and irreducible curve $C$ with $h^1(X, \C C)=1$, and
either $r=1$ or $\C C^2=0$. [This follows
if the linear system of sections of 
\C H is composed with a pencil by \ir{Bertini}(a). If not,
then $\C H^2>0$, $r=1$ by \ir{Bertini}(b), and
$\C H=\C C$ is irregular by \ir{siposirrzero}, so
$h^1(X, \C H)=h^1(X, \C C)=1$ by (a).] Thus
$h^0(X, \C H)=r^2C^2/2+r+1$ by Riemann-Roch. 

By \ir{keycor}, the components of $N$ are smooth and rational,
of negative self-intersection, and each component $M$ which is
not also a fixed component of the linear system of sections of $-K_X$ has by adjunction
either $M^2=M\cdot K_X=-1$ or $M^2=-2$ and $M\cdot K_X=0$. 
Moreover,
each connected component of $N$ meets \C H
positively, by \ir{RFb}. Note that an irreducible component
of $N$ having positive intersection product with \C H cannot be a 
component of $D$, since $-K_X\cdot \C H=0$ means
$D$ is disjoint from a general section of \C H. 
In particular, each connected component of $N$ 
is either disjoint
from $D$ (and hence each irreducible component
of such a connected component has self-intersection $-2$) 
or contains a string of components 
$N_1,\ldots,N_t$, none being 
a component of $D$, with (we may assume) $N_i\cdot N_{i+1}>0$,
$C\cdot N_1>0$, $-K_X\cdot N_i=0$ for $i<t$, and $-K_X\cdot N_t=1$.

Consider an irreducible component $M$ of $N$ which meets \C H,
and hence which is not a 
component of $D$. Clearly, $M$ is a 
fixed component of $|rC+M|$ (so $rC+M$ and $rC$ move in complete 
linear systems of the same dimension). Also the general section of 
$\C O_X(rC+M)$ is connected (since $C\cdot M>0$),
reduced, and has no fixed components in common with $D$.
Thus $h^1(X, \C O_X(rC+M))=0$ by \ir{RFa}(a) if $M^2=-1$,  
and $h^1(X, \C O_X(rC+M))=1$ by \ir{RFa}(b) if $M^2=-2$.
We can in either case now compute $h^0(X, \C O_X(rC+M))$ by
plugging into Riemann-Roch; setting the result equal 
to $h^0(X, \C H)=r^2C^2/2+r+1$ and solving gives $C\cdot M=1$. 

We now consider two cases. First say 
the connected components of $N$ are not all disjoint
from $D$; thus there is a string 
$N_1,\ldots,N_t$, as above, with $N_1$ meeting \C H and 
$N_t$ meeting $D$. Note that
$L=rC+N_1+\cdots+N_t$ is numerically effective, since it meets
each $N_i$ nonnegatively, and regular by \ir{RFa}(a). If $N_1+\cdots+N_t$ is not all of \C N, 
then the difference 
$N'$ is nontrivial. Since $h^0(X, \C O_X(L))=h^0(X, \C F)$, we see from
$0\to \C O_X(L)\to \C F\to \C F\otimes\C O_{N'}\to 0$ that
$h^0(N', \C F)=0$, contradicting \ir{keycor}. I.e., $N'$ must be trivial,
\C N is precisely the class of $N_1+\cdots+N_t$, and \C F is the class of $L$,
hence regular. If $t=1$, we are done with this case. If $t>1$, to see that
$C\cdot N_i=0$ for $i>1$, that $N_i\cdot N_j=0$ for $j>i+1$, and that $N_i\cdot N_{i+1}=1$
for each $i$,
note by \ir{RFa} that $h^1(X, \C O_X(rC+N_1+N_2))$ is 1 if $N_2^2=-2$        
and 0 if $N_2^2=-1$. But $h^0(X, \C O_X(rC+N_1+N_2))=h^0(X,\C H)=r^2C^2/2+r+1$.
Using Riemann-Roch and our knowledge of $h^1$ to compute
$h^0(X, \C O_X(rC+N_1+N_2))$ explicitly by setting it equal to
$r^2C^2/2+r+1$ and simplifying, 
gives $(rC+N_1)\cdot N_2=1$. Since $N_1\cdot N_2>0$, we see
$N_1\cdot N_2=1$ and $C\cdot N_2=0$. This argument can be repeated with
$rC+N_1+\cdots+N_i$ for each $i\le t$, to give the result.

Consider now the case that the connected components of $N$ are all disjoint
from $D$. This means that $\C F\otimes\C O_D=\C O_D$,
and so from $0\to \C F+K_X\to \C F\to \C O_D\to 0$ 
and \ir{ODisirr} we see \C F is irregular. 
Let $M$ be an irreducible
component of $N$ meeting \C H. Then $M^2=-2$, and 
$h^1(X, \C O_X(rC+M))=1$ (as observed above). 
If $M$ is not all of \C N, the difference 
$N'$ is nontrivial. Thus we can consider
$0\to \C O_X(rC+M)\to \C F\to \C F\otimes\C O_{N'}\to 0$. 
Since \C F is irregular and $h^1(N', \C F)=0$ (\ir{keycor}), 
the $1$-dimensional space $H^1(X, \C O_X(rC+M))$ maps
in the long exact sequence
isomorphically to $H^1(X, \C F)$.
Thus $h^0(N', \C F)=0$, contradicting \ir{keycor}. I.e., $N'$ must be trivial 
and \C N is precisely the class of $M$. Thus 
\C F is the class of $rC+M$, hence $h^1(X, \C F)=1$, and
$r>1$ (since otherwise $\C F\cdot M<0$). \qed

The following result is needed to deal with cases in which the base locus contains
an anticanonical divisor.

\prclm{Lemma}{nefadjoint}{Let \C F be a numerically effective class on a
smooth projective rational anticanonical surface $X$.
Let \C N be the class of the fixed part of the linear system of sections of \C F.
\itemitem{(a)} There is a birational morphism $X\to Y$ to a smooth surface $Y$
such that \C F is the pullback
of a numerically effective class $\C L$ on $Y$, and $\C L\cdot \C E>0$ 
for every class
\C E of an irreducible exceptional divisor on $Y$.
\itemitem{(b)} Let $X\to Y$ be a birational morphism to a smooth surface $Y$ 
such that \C F is the pullback of a class $\C L$ on $Y$.
Let $\C M$ be the class of the fixed part of the linear system of sections of $\C L$. 
Then $h^0(X, \C N+K_X)>0$ if and only if $h^0(Y, \C M+K_Y)>0$.
\itemitem{(c)} If  $h^0(X, \C F+K_X)>0$ and if
$\C F\cdot \C E>0$ for every class \C E of an irreducible
exceptional curve, then $\C F+K_X$ is numerically effective.}

\Prf (a) If \C E is the class of an 
irreducible exceptional divisor $E$ with
$\C F\cdot E=0$, let $X\to Y'$ be the birational morphism contracting
$E$. Then the corresponding homomorphism $\hbox{Pic}(Y')\to \hbox{Pic}(X)$
is an isomorphism of Pic$(Y')$ to $E^\perp$. Thus there is a class 
$\C F'$ on $Y'$ whose pullback to $X$ is \C F. 
Note that $\C F'$ is numerically effective by \ir{one}.
We can repeat this process unless $Y'$ has no exceptional classes
perpendicular to \C F, which must eventually occur since
the rank of Pic$(Y')$ is less than the rank of Pic$(X)$.

(b) Since any birational morphism of surfaces factors into a sequence
of contractions of irreducible exceptional divisors [Ha, V.3.2], 
it is enough by induction to show  $h^0(X, \C N+K_X)>0$ if 
and only if $h^0(Y, \C M+K_Y)>0$ in the case that
$X\to Y$ is obtained by contracting a single irreducible exceptional 
curve $E$,
whose class we denote \C E. Thus we have $\C F\cdot \C E=0$. Also,
the pullback of $K_Y$ is $K_X-\C E$ [Ha, V.3.3]. Thus, by \ir{one},
$h^0(X, \C F+K_X-\C E)=h^0(Y, \C L+K_Y)$ and 
$h^0(Y, \C L)=h^0(X, \C F)$. 

Say $h^0(X, \C N+K_X)>0$. Since $\C F-\C N\in\hbox{EFF}$, so is $\C F+K_X$.
But $(\C F+K_X)\cdot \C E=-1$ so $E$ is in the fixed part of the linear system of sections of
$\C F+K_X$. Of course, a section of $-K_X$ is in the fixed part of the linear system of sections of \C F.
Thus $h^0(X, \C F)=h^0(X, \C F+K_X)=h^0(X, \C F+K_X-\C E)$ and hence
$h^0(Y, \C L)=h^0(Y, \C L+K_Y)$, which means that a section of $-K_Y$ is in 
the fixed part of the linear system of sections of \C L; i.e., $h^0(Y, \C M+K_Y)>0$. 

Conversely, say $h^0(Y, \C M+K_Y)>0$. Then
$h^0(Y, \C L)=h^0(Y, \C L+K_Y)$, so 
$h^0(X, \C F)=h^0(X, \C F+K_X-\C E)$. Thus a section of $-K_X+\C E$ is
in the fixed part of the linear system of sections of \C F so certainly
a section of $-K_X$ is in the fixed part, whence $h^0(X, \C N+K_X)>0$.

(c) Let $C$ be a reduced and irreducible curve on $X$.
If $C^2\ge 0$ or $C\cdot K_X\ge 0$, then clearly $C\cdot(\C F+K_X)\ge 0$.
So say $C^2<0$ and $C\cdot K_X<0$. By adjunction, $C$ is an exceptional curve;
i.e., a smooth rational curve with $C^2=C\cdot K_X=-1$. By assumption,
$\C F\cdot C>0$, so again $C\cdot(\C F+K_X)\ge 0$. \qed

By the following lemma we see that there are two ways 
a numerically effective divisor can
contain an anticanonical divisor in its 
fixed part. One of these is described explicitly; the other
case reduces (essentially by subtracting off
the anticanonical divisor) to cases previously worked out, in which
the fixed part does not contain an anticanonical divisor. An explicit
description could be made in this second case, too. We chose not to include
it here, since the statement
of the lemma is already somewhat complicated and since it
is not hard using \ir{mainthm} to work out an explicit description
if one chooses.

\prclm{Lemma}{DinN}{Let \C F be a numerically effective class on a
smooth projective rational anticanonical surface $X$,
let $D$ be a nonzero section of $-K_X$,
let \C N denote the class of the fixed part of the linear system of sections of \C F and let $\C H=\C F-\C N$.
Then $\C N+K_X\in\hbox{EFF}$ if and only if $\C F\cdot D=0$ but
$\C F\otimes \C O_D$ is nontrivial. In this case, there exists a birational
morphism of $X$ to a smooth
projective rational anticanonical surface $Y$, 
and either: $K_Y^2<0$, there is a
numerically effective class $\C F'$ on $Y$, \C F is the pullback of $\C F'-K_Y$, and
$0=h^1(Y, \C F')=h^1(X, \C F)$; 
or $K_Y^2=0$, \C H and \C N are the pullbacks
of $-sK_Y$ and $-rK_Y$ for some integers $s\ge 0$ and $r>0$ respectively, 
and $h^1(X, \C F)=\sigma$, where $\sigma=0$ if 
$s=0$, and otherwise $r<\tau$ and $\sigma=s/\tau$, where
$\tau$ is the least positive integer such that 
the restriction of $-\tau K_X$ to $D$ is trivial.}

\Prf Suppose a section of $-K_X$ occurs in the 
fixed part of the linear system of sections of \C F; i.e., $h^0(X, \C N+K_X)>0$. In particular,
this means that $h^0(X, -K_X)=1$, and thus that $(-K_X)^2\le 0$ 
(otherwise,
by Riemann-Roch, $h^0(X, -K_X)>1$). 

We wish to prove that $\C F\cdot K_X=0$. Denote $X$ by $X_1$,
$\C F$ by $\C F_1$ and \C N by $\C N_1$. By \ir{nefadjoint}(a), 
there is a birational morphism $X\to Y$ such that \C F is the pullback
of a class $\C L_1$ (numerically effective by \ir{one}(d))
on $Y$ which meets every irreducible exceptional class on $Y$
positively; by \ir{nefadjoint}(b), $h^0(X, \C M_1+K_Y)>0$, where 
$M_1$ is the fixed part of the linear system of sections of $\C L_1$. By \ir{nefadjoint}(c), $\C L_1+K_Y$
is numerically effective. Denote $Y$ by $X_2$,
$\C L_1+K_Y$ by $\C F_2$ and the class of the fixed part of 
the linear system of sections of $\C F_2$ by $\C N_2$. 

Thus given a numerically effective $\C F_i$ on $X_i$ such that 
$h^0(X_i,\C N_i+K_{X_i})>0$, where $\C N_i$ is the class of the
fixed part of the linear system of sections of $\C F_i$, we obtain a numerically
effective class $\C L_i$ on a surface $X_{i+1}$ such that
$h^0(X_{i+1},\C M_i+K_{X_{i+1}})>0$, where $\C M_i$ is the class of the
fixed part of the linear system of sections of $\C L_i$ and 
every irreducible exceptional class on $X_{i+1}$ meets $\C L_i$ positively, 
and from $\C L_i$ we obtain
a numerically effective class $\C F_{i+1}$ on $X_{i+1}$. 
Since $-K_{X_i}\in\hbox{EFF}$ for each $i$, it meets $\C F_i$ nonnegatively. 
Since this process of going from $i$ to $i+1$ 
reduces the number of integral divisors in the fixed part, it eventually
must stop. I.e., for some $j$, $\C F_j$ is numerically effective
but $h^0(X_j, \C N_j+K_{X_j})=0$. Since for each $i\le j$ we have
$h^0(X_i, \C M_{i-1}+K_{X_i})>0$ (which, as remarked above, is 
impossible if $K_{X_i}^2>0$), we see that 
$K_{X_i}^2\le 0$ for all $i\le j$. 

Denote $K_{X_i}$ by $K_i$. Then for $i<j$,
$\C F_i\cdot K_i=\C L_i\cdot K_{i+1}$, and
$\C F_{i+1}=\C L_i+K_{i+1}$
so $\C F_{i+1}\cdot K_{i+1}=\C F_i\cdot K_i+K_{i+1}^2$.
But $0\ge K_{i+1}^2\ge K_i^2$ (with $K_{i+1}^2=K_i^2$ precisely if
$X_i=X_{i+1}$, which would just mean that $\C F_i=\C L_i$ already 
meets every irreducible exceptional class on $X_i$
positively); thus $0\ge \C F\cdot K_X>\C F_j\cdot K_j$ unless
$X_2=X_j$ and $K_2^2=0$, in which case $\C F_j=\C L_1+(j-1)K_2$.
We see it is enough to consider two cases: either 
$0>\C F_j\cdot K_j$, or $0=\C F_j\cdot K_j$. In the latter case,
$X_2=X_j$, $K_2^2=0$ and by \ir{Lattice}, $\C F_j=-sK_2$ for some
$s\ge 0$, so for $r=j-1$, $\C L_1=-sK_2-rK_2$.

Consider first the contingency  $0>\C F_j\cdot K_j$.
Let $D_j$ be a nontrivial section of $-K_j$ (e.g., take the image of 
$D$ under $X_1\to X_j$). 
Consider $0\to \C F_j\to \C L_{j-1}\to \C L_{j-1}\otimes\C O_{D_j}\to 0$.
From \ir{sizero} and \ir{Bertini}(a), \ir{siposirrzero}, and \ir{siposhirr},  
we see that a numerically effective divisor not containing an anticanonical divisor in the
fixed part of its linear system of sections but meeting the anticanonical class positively is regular;
i.e., $\C F_j$ is regular, so the sequence is exact on global
sections and on $h^1$. Since $\C F_j$ and $\C L_{j-1}$ differ only
in the fixed components of their linear systems of sections, they have isomorphic $H^0$'s.
From Riemann-Roch, $h^0(X_j, \C F_j)=(\C F_j^2-K_j\cdot\C F_j)/2+1$
and $h^0(X_j, \C L_{j-1})=
(\C L_{j-1}^2-K_j\cdot\C L_{j-1})/2+1+h^1(X_j, \C L_{j-1})$.
Setting these equal and simplifying using $\C F_j=\C L_{j-1}+K_j$
gives $\C L_{j-1}\cdot K_j=h^1(X_j, \C L_{j-1})$, but
$\C L_{j-1}\cdot K_j\le 0$ since $\C L_{j-1}$ is numerically
effective, so $h^1$ vanishes for $\C L_{j-1}$.
Thus $(\C F_j-K_j)\cdot K_j=\C L_{j-1}\cdot K_j=0$, or $\C F_j\cdot 
K_j=K_j^2$.  Since 
$\C F\cdot K_1+K_2^2=\C F_2\cdot K_2\ge \C F_j\cdot K_j=K_j^2$, 
we have $\C F\cdot K_X=\C F\cdot K_1\ge  K_j^2-K_2^2\ge 0$. 
Since \C F is numerically effective,
we see that $\C F\cdot K_X=0$ as desired. 
We have also proved that in this case
$K_j^2=K_2^2$, hence that $X_j=X_2$, so $\C F_j=\C L_1+(j-1)K_2$. 
Substituting into $\C F_j\cdot K_j=K_j^2$ gives 
$(\C L_1+(j-1)K_2)\cdot K_2=K_2^2$ or 
$\C L_1\cdot K_2=(2-j)K_2^2$ 
and so we have
$0=\C F\cdot K_X=\C L_1\cdot K_2=(2-j)K_2^2$; i.e.,
either $j=2$ or $K_2^2=0$. But 
$K_2^2=0$ implies $\C F_j\cdot K_j=0$,
contradicting the hypothesis $0>\C F_j\cdot K_j$.
Thus $j=2$, $K_2^2<0$, and $Y$ in the statement
of \ir{DinN} is $X_2$, while $\C F'$ is $\C F_2$.

We note that in the course of the argument above 
we found that $\C L_{j-1}=\C L_1$ is 
regular; thus \C F, being a pullback
of $\C L_1$, is regular too. But $h^2(X, \C F+K_X)=0$ by duality,
so $h^1(D, \C F)=0$ follows from regularity of \C F and from
$0\to \C F+K_X\to \C F\to\C F\otimes\C O_D\to 0$. From
\ir{ODisirr} we therefore see that $\C F\otimes\C O_D$ is nontrivial.

Now consider the remaining contingency:
$X_2=X_j$, $K_2^2=0$, $\C F_j=-sK_2$ where $s\ge 0$, 
and $\C L_1=-sK_2-rK_2$ for $r=j-1$. Since
$\C F_j$ is numerically effective of self-intersection $0$
and the fixed part of its linear system of sections
does not contain a section of $-K_2$, 
either \ir{sizero} applies or $\C F_j=0$. Either way, $\C N_j=0$,
with $h^1(X_2, \C F_j)=\sigma$,
where $\sigma=0$ if $\C F_j=0$ and otherwise by \ir{Bertini}(a) 
$\C F_j=\sigma\C C$ for some positive
integer $\sigma$, where \C C is the class of some reduced and 
irreducible curve $C$ moving in a pencil. 

Clearly, if $s=0$, then $\sigma=0$ and $\C L_1=-rK_2$; i.e., we can take $Y=X_2$,
and then \C H is the pullback of $-sK_2$ with $s=0$, and \C N is the pullback
of $-rK_2$, as the statement of \ir{DinN} requires. Moreover, in this case 
$h^0(X_2, \C L_1)=1$, so $h^1(X_2, \C L_1)=0$ by Riemann-Roch, so 
$h^1(X, \C F)=h^1(X_2, \C L_1)=\sigma=0$, as claimed. 

If $s>0$, then $\C C=-\tau K_2$, where $\sigma\tau=s$. Thus
the sections of $\C F_j=\sigma\C C=-sK_2$ are base point free and we have
$h^1(X_2, \C F_j)=\sigma$. But $h^0(X_2, \C F_j)=h^0(X_2, \C L_1)$
and by Riemann-Roch $h^0-h^1$ is the same for $\C F_j$ and $\C L_1$,
and hence $h^1$ is the same too, so
$h^1(X_2, \C L_1)=h^1(X, \C F_j)=\sigma$. Moreover, since the pullback 
of $\C F_j$ to $X$ has the same $h^0$ as \C F, we see that
\C H is the pullback of $\C F_j=-sK_2$ and hence \C N is the pullback
of $-rK_2$, and so again $Y=X_2$.

We also have $r<\tau$, since 
$-rK_2$ is the class of a unique effective divisor but the sections of $\C C=-\tau K_2$ move in a pencil.
Let $t<\tau$; then $\C C+tK_2=-(\tau-t)K_2\in\hbox{EFF}$.
Since the sections of $\C C=-\tau K_2$ are fixed component free and hence 
$D_2\cdot\C C=0$,
the restriction of $-K_2$ to $C$ is trivial, so
from $0\to  -\C C-tK_2\to -tK_2\to -tK_2\otimes\C O_C\to 0$
we see $h^0(X_2, -tK_2)=1$ and, by Riemann-Roch, $h^1(X_2, -tK_2)=0$. 
But, denoting the pullback of $K_2$ to $X$ also by $K_2$,
we have $0=h^1(X_2, -tK_2)=h^1(X, -tK_2)$, so from 
$0\to  -tK_2+K_X\to -tK_2\to -tK_2\otimes\C O_D\to 0$, we see 
$h^1(D, -tK_2)=0$; i.e., for $t<\tau$, the restriction of $-K_2$ to 
$D$ is, by \ir{ODisirr}, nontrivial. But $C$ is disjoint from $D_2$, hence
the pullback of \C C has trivial restriction to $D$; thus $\tau$ is the least positive 
multiple of $-K_2$ whose restriction to $D$ is trivial. 
Moreover, it follows that
the restriction of \C H to $D$ is trivial, and hence that
$\C F\otimes\C O_D=-rK_2\otimes\C O_D$ is nontrivial.

Conversely, say $\C F\cdot D=0$. 
By \ir{sizero}, \ir{siposirrzero}, and \ir{siposhirr}, if $\C N+K_X\notin\hbox{EFF}$ 
(i.e., if the fixed part of the linear system of sections of \C F does not contain 
an anticanonical divisor),
then no fixed component of the linear system of sections of \C F is a fixed 
component of the linear system of sections of $-K_X$. In this case
$\C F\cdot D=0$ means that \C F has a section which is disjoint from
$D$, and hence that $\C F\otimes \C O_D$ is indeed trivial. Thus, if
$\C F\otimes \C O_D$ is nontrivial, it must be that
$\C N+K_X\in\hbox{EFF}$.  \qed

Our results on base points depend on the following lemma.

\prclm{Lemma}{folklore}{Let $C$ be an integral projective curve whose
dualizing sheaf $K_C$ is locally free rank $1$. If $h^0(C, K_C)>0$
(i.e., $C$ is not smooth and rational), then $K_C$ is generated by global sections.}

\Prf  This is a special case of Theorem D of [Cn]. \qed

The proof in [F] corresponding to the next result assumes that
$-K_X$ has a reduced section, but this can be avoided.
 
\prclm{Lemma}{basepoints}{Let \C F be the class of an effective fixed 
component free divisor on a
smooth projective rational anticanonical surface $X$. Let $D$ be a 
nonzero section of
$-K_X$. Then the sections of \C F have a base point if and only if $\C F\cdot D=1$,
in which case the base point is unique and lies on $D$.}

\Prf First we check that there is no base point if
$\C F\cdot D\ne 1$. Since \C F is numerically effective,
we have $\C F \cdot D\ge 0$. Say $\C F \cdot D=0$; if also
$\C F^2=0$ then, being fixed component free, the sections of \C F are also
obviously base point free (since distinct sections of \C F are disjoint).
So say $\C F^2>0$. A general section $F$ of \C F is reduced
and irreducible by \ir{Bertini}. Thus we may assume 
$F$ and $D$ are disjoint, and
so $\C O_F\otimes \C F=K_F$ by adjunction. From
$0\to \C O_X\to \C F\to \C F\otimes\C O_F\to 0$, we see 
the sections of \C F surject to those of $\C F\otimes\C O_F=K_F$. 
By \ir{folklore}, $K_F$ is generated by global sections (and hence
the sections of \C F are base point free) unless $F$ is rational; i.e., unless
$F^2=-2$, contrary to assumption.

Now say $\C F \cdot D>1$. Thus $\C O_F\otimes \C F=-K_X\otimes K_F$
has degree at least $2g$, where $g$ is the genus of $F$. By 
[D, Proposition 7, p.\ 59], $\C O_F\otimes \C F$ is generated by
global sections, and hence as before the sections of \C F are base point free.

Conversely, suppose $\C F\cdot D=1$. Then $\C F^2$ is odd by adjunction,
so positive, so a general section $F$ of \C F is integral by \ir{Bertini},
and we may assume $F$ and $D$ have no common components and hence meet at a
single point, $x$, which must be smooth on both $F$ and $D$. 
Let $Y\to X$ be the morphism resulting from
blowing up $x$ and let \C E be the class of the exceptional
locus $E$. Then $\C F'=\C F-\C E$ is numerically effective
and its sections fixed component free on $Y$, and the proper transform $D'$ of 
$D$ is anticanonical, but $\C F'\cdot D'=0$. Thus, by our foregoing argument,
the sections of $\C F'$ are base point free. By \ir{Bertini}(a) (if ${\C F'}^2=0$) or 
\ir{siposirrzero} and \ir{siposhirr}(a) (if ${\C F'}^2>0$), we see $h^1(Y, \C F')=1$. 
But $\C F'+\C E$ is numerically effective (since $\C F'$ is and 
since it meets \C E nonnegatively) of positive self-intersection
and meets the anticanonical class,
so from \ir{siposirrzero} and \ir{siposhirr}, we see it must be regular.
Thus from $0\to \C F'\to \C F'+\C E\to (\C F'+\C E)\otimes\C O_E\to 0$
we see that $h^0(Y, \C F')=h^0(Y, \C F'+\C E)$, so \C E is the class of a fixed
component of the linear system of sections of $\C F'+\C E$. Since $\C F'+\C E$ is the pullback
of \C F to $Y$, we see that the sections of \C F have a base point at $x$. \qed

\noindent{\it Proof} of \ir{mainthm}: (a) By \ir{DinN}, we see that the fixed part
of the linear system of sections of \C F cannot contain 
an anticanonical divisor. By \ir{sizero}, \ir{Bertini}, \ir{siposirrzero} and
\ir{siposhirr}, we see that \C F is regular and its linear system of sections
is fixed component free, and so
base point free by \ir{basepoints}.

(b) Here \C F must have positive self-intersection (to satisfy adjunction),
so, by \ir{siposirrzero} and \ir{siposhirr}, \C F is regular. If its linear system
of sections is also
fixed component free, then by \ir{basepoints} it has a unique base point, 
on $D$. If the sections of \C F have a fixed component, then the result follows by
\ir{siposirrzero} and \ir{siposhirr}, again.         

Conversely, say $\C F=\C H +\C N$, with $\C H=r\C C$, $h^1(X, \C H)=r$ and 
$\C N=\C N_1+\cdots+\C N_t$, as in the statement of (b).
It is easy to check that $\C N^2=\C N\cdot K_X=-1$. Since \C F is regular, 
we have using $\C F=\C H+\C N$, $\C C\cdot \C N=1$ and Riemann-Roch that 
$h^0(X, \C F)= (\C H^2-\C H\cdot K_X)/2+r+1$. But $h^1(X, \C H)=r$ by hypothesis
so we now see $h^0(X, \C F)=h^0(X, \C H)$. Thus $\C N$ is the class of the fixed part of
the linear system of sections of $\C H+\C N$, as we needed to show.

(c) If the sections of \C F are fixed component free, then first they are base point free
by \ir{basepoints}, and second either 
\ir{siposhirr}(a) applies and gives the result or 
\ir{sizero} applies and, with \ir{Bertini}, gives the result. 
If the sections of \C F are not fixed 
component free but the fixed part does not contain an anticanonical divisor, then only
\ir{siposhirr}(b) applies, giving the result. 

(d) This is just \ir{DinN}.  \tempqed

\noindent {\it Proof} of \ir{allornone}: Let \C F be a numerically effective class
such that the fixed part of its linear system of sections is nontrivial but 
does not contain an anticanonical divisor. 
Of \ir{sizero}, \ir{siposirrzero} and \ir{siposhirr}, only \ir{siposhirr}(b)
applies, whence no fixed component of the linear system of sections of 
\C F is a component of an arbitrary anticanonical divisor $D$.
\tempqed

We now prove \ir{intr}, of the Introduction:
\vskip\baselineskip

\Prf So \C F is a numerically effective divisor class on
a smooth projective rational surface $X$ with
an effective anticanonical divisor $D$. Then $h^0(X, \C F)>0$
by \ir{NefIsEffCor}. Moreover, by \ir{regcond}(a), $h^1(X, \C F)>0$ if and only if
$\C F\cdot K_X=0$ and a general section of $\C F$ has a connected 
component disjoint from $D$. Thus we only need to show that 
$1+h^1(X, \C F)$ is the number of connected
components of a general section of $\C F-K_X$.

By \ir{ODisirr}, anticanonical divisors are connected, hence 
the result follows for $\C F=\C O_X$, so we may assume that
\C F is nontrivial. Suppose that $h^1(X, \C F)=0$. Then
$h^1(X, -(\C F-K_X))=0$ by duality,
so from $0\to -(\C F-K_X)\to \C O_X\to \C O_C\to 0$,
where $C$ is a general section of $\C F-K_X$, we
see $C$ is connected, as desired. Suppose $h^1(X, \C F)>0$;
then, by \ir{mainthm} and \ir{allornone}, $-K_X\cdot\C F=0$  
and if the general section of \C F is not disjoint from $D$ then,
as we shall first assume,
there is a birational morphism $X\to Y$ to a smooth surface $Y$ with
$K_Y^2=0$ and (regarding pullbacks from $Y$ as classes on $X$)
$\C F=-\tau\sigma K_Y-rK_Y$, where $1\le\sigma$, $1\le r<\tau$, 
$-rK_Y$ is the class of the fixed part of the linear system of sections of \C F, and $-\tau K_Y$
is a pencil with which the linear system of sections of $-\tau\sigma K_Y$ is composed.
Suppose $X=Y$; if $r+1<\tau$ then $\C F-K_X=-\tau\sigma K_Y-(r+1)K_Y$,
and $h^0(Y,-(r+1)K_Y)=1$, so a general section of $\C F-K_X$
has $1+\sigma$ connected components, while $h^1(X, \C F)=\sigma$
by \ir{mainthm}(c), as desired. If $r+1=\tau$, then 
the linear system of sections of
$\C F-K_X=-\tau(\sigma+1)K_Y$
is fixed component free, and a general section again has 
$1+\sigma$ connected components, as desired. If $X\to Y$ is not the identity,
then $X$ is a blowing up of points (possibly infinitely near)
of $Y$. In this case, $\C F-K_X=-\tau\sigma K_Y-rK_Y-K_X$, 
hence a general section once more has
$1+\sigma$ connected components, as needed, if we check that $h^0$ of
$-rK_Y-K_X$ is 1. Since in any case $h^0$ of $-rK_Y$ is 1,
we have $r<\tau$; thus $h^0$ of $-rK_Y-K_X$ is 1 if $r+1<\tau$, since $h^0$ of
$-rK_Y-K_X+(-K_Y+K_X)=-rK_Y-K_Y$ is 1 and $-K_Y+K_X$, being
a sum of classes of exceptional curves, is the class of an effective divisor. If $r+1=\tau$,
then the linear system of sections of $-rK_Y-K_Y$ is a pencil, 
and adding $-(-K_Y+K_X)$, which is the class of an antieffective divisor, we get 
the class $-rK_Y-K_X$, which has a unique effective section.
 
So now we may assume that the general section of \C F is
disjoint from $D$.
If the linear system of sections of \C F is fixed component free, then, by \ir{mainthm}(c) 
and \ir{Bertini}, $h^1(X,\C F)$ equals the number of components
of a general section $F$ of \C F, and $F$ is disjoint
from $D$. Thus $\C O_F\otimes(\C F-K_X)=\C O_F\otimes\C F$,
so from $0\to \C O_X\to \C F\to \C O_F\otimes\C F\to 0$
and $0\to -K_X\to \C F-K_X\to \C O_F\otimes(\C F-K_X)\to 0$
we see that the sections of $\C F-K_X$ are just sums
of sections of \C F and $-K_X$, hence the number of components
of a general section of \C F is one less than the number of components
of a general section of $\C F-K_X$, as desired. If the fixed part $N$
of the linear system of sections of \C F is not trivial but does not contain
an anticanonical divisor,
then by \ir{mainthm}(c) $h^1(X, \C F)=1$ and the general section
of \C F is connected but, by \ir{allornone}, disjoint from $D$. 
Arguing as in the previous case
shows that the general section of $\C F-K_X$ has two
connected components, as required. \qed

\rem{Remark}{remfortony} Our results have mostly concerned 
numerically effective classes, but one may be interested in 
arbitrary classes. In general, to determine
dimensions of complete linear systems, fixed components and
base points, 
given arbitrary divisor classes on a smooth
projective rational anticanonical surface, one essentially needs to know
the monoid EFF of effective classes on $X$.

Since EFF determines the cone NEF of numerically effective classes, given
an arbitrary class \C F and knowing EFF one can determine if $\C F\in\hbox{EFF}$
and, if so, find a decomposition $\C F=\C L+\C M$, where \C L is maximal with respect to
$\C L\in \hbox{NEF}$ and $\C L\le \C F$ (where $\C L\le \C F$ means $\C F-\C L\in\hbox{EFF}$).
Then $h^0(X, \C F)=h^0(X,\C L)$, and the results obtained herein for
numerically effective classes apply to \C L.

To make this discussion concrete, consider the case that $X$ is a blowing up of \pr 2.
(Blowings up of other relatively minimal models can with minor changes be handled
similarly.) In this case the following data suffice to determine EFF:
an {\it exceptional configuration\/} $\{\C E_0,\ldots,\C E_n\}\subset \hbox{Pic}(X)$ 
(i.e., a basis  $\C E=\{\C E_0,\ldots, \C E_n\}$ of
Pic$(X)$, where $\C E_0$ is the pullback of the class of a line with respect to
some birational morphism $X\to \pr 2$
and the other classes $\C E_i$ are the classes of the exceptional loci
corresponding to a factorization of $X\to \pr 2$ into a sequence of 
monoidal transformations); the kernel $\Lambda\subset \hbox{Pic}(X)$ of the functorial
homomorphism $\hbox{Pic}(X)\to \hbox{Pic}(D)$, where $D$ is an effective anticanonical divisor;
and the fixed components of $|D|$.

In fact, EFF is generated by: (a) the fixed components of $|D|$; (b) classes \C F
such that $\C F^2=\C F\cdot K_X=-1$ and $\C F\cdot \C E_0\ge 0$; (c) classes \C F
such that $\C F=\C E_i-\C E_j$ where $0<i<j$ and $\C F\in \Lambda$; (d) classes \C F
such that $\C F^2=-2$, $\C F\cdot K_X=0$, $\C F\cdot \C E_0>0$ and $\C F\in \Lambda$;
and (e) classes \C F such that $\C F^2-\C F\cdot K_X\ge 0$ and $\C F\cdot \C E_0\ge 0$.

To see this, we first check that classes of each of these types are in EFF.
This is obvious for classes of type (a). For a class \C F of one of the remaining types
note that $h^2(X,\C F)=h^0(X, K_X-\C F)=0$ since $\C E_0$ is numerically effective
but $(K_X-\C F)\cdot \C E_0<0$. Now effectivity follows in cases (b) and (e) by \ir{RR}.
In cases (c) and (d), by \ir{RR} we have $h^0(X, \C F)=h^1(X, \C F)$; we also have
$h^2(X, K_X+\C F)=0$: in case (d) since $h^2(X, K_X+\C F)=h^0(X, -\C F)$ and 
$-\C F\cdot\C E_0<0$, and in case (c) since, if $h^2(X, K_X+\C F)>0$, then 
$-\C F\in\hbox{EFF}$, which is impossible since $0<i<j$ implies $\C E_j-\C E_i$
cannot be in EFF (as exceptional divisors, $E_j$ may be contained in $E_i$
but not vice versa). In either case (c) or (d), since $\C F\in \Lambda$, the restriction
of \C F to $D$ is trivial so
we have $0\to \C F+K_X\to \C F\to \C O_D\to 0$, hence
$h^0(X, \C F)=h^1(X, \C F)\ge h^1(D, \C O_D)=1$.

Conversely, if \C F is the class of a reduced irreducible curve of negative self-intersection,
then $\C F\cdot \C E_0\ge 0$ and either \C F 
is a fixed component of $|D|$ (which is case (a)), or \C F is 
not the class of a fixed component of $|D|$ and so meets $-K_X$ nonnegatively. Hence by adjunction either
$\C F^2=\C F\cdot K_X=-1$ (case (b)), or $\C F^2=-2$, $\C F\cdot K_X=0$ and \C F is
disjoint from $D$, and thus in $\Lambda$
(which is case (c) if $\C F\cdot \C E_0=0$, and case (d) if $\C F\cdot \C E_0>0$).
In general, if \C F is in EFF, we can write $\C F=\C L+\C N$, where \C N
is a nonnegative sum of curves of negative self-intersection (and hence in the monoid
generated by types (a) through (d)), and \C L is numerically effective
(and hence of type (e)).\endremark

\Refs
\widestnumber\key{Ma}

\bibitem{A} M. Artin\paper Some numerical criteria 
for contractability of curves
on algebraic surfaces
\jour Amer.\ J.\  Math.
\vol 84 
\yr 1962
\pages 485--497
\endref

\bibitem{Cn} F.\ Catanese \paper Pluricanonical Gorenstein curves \paperinfo preprint
\endref

\bibitem{Co} C. Ciliberto \paper On the degree and genus of smooth
curves in a projective space
\jour Adv.\ Math. 
\vol 81 
\yr 1990
\pages 198--248
\endref

\bibitem{D} M. Demazure \paper Surfaces de Del Pezzo - IV. Syst\`emes anticanoniques
\inbook S\'eminaire sur les Singularit\'es des Surfaces
\vol 777 \bookinfo LNM
\yr 1980
\endref

\bibitem{F} R. Friedman \paper Appendix: Linear systems on anticanonical pairs
\pages 162--171 \inbook The Birational Geometry of Degenerations
\vol Prog.\ Math.\ 29 \publ Birkhauser \publaddr Boston \yr 1983
\endref

\bibitem{\trans} B. Harbourne
\paper Complete linear systems on rational surfaces
\jour Trans.\ A.\ M.\ S.\ 
\vol 289 
\yr 1985
\pages 213--226
\endref

\bibitem{\annalen} B. Harbourne \paper Very ample divisors on rational surfaces
\jour Math.\ Ann. \vol 272 \yr 1985 \pages 139--153
\endref

\bibitem{\duke} B. Harbourne
\paper Blowings-up of \pr 2 and their blowings-down
\jour Duke Math. J. 
\vol 52 
\yr 1985
\pages 129--148
\endref

\bibitem{\bowdoin} B. Harbourne \paper Automorphisms of $K3$-like rational surfaces
\inbook Proc.\ Symp.\ Pure Math. \vol 46
\yr 1987 \pages 17--28
\endref

\bibitem{\sundancetwo} B. Harbourne \paper Automorphisms of cuspidal K3-like surfaces
\inbook Algebraic Geometry: Sundance 1988 \pages 47--60 \vol 116
\bookinfo Contemporary Mathematics
\yr 1991
\endref

\bibitem{\mtnwest} B. Harbourne \paper Rational surfaces with $K^2>0$
\paperinfo to appear, Proc. A.M.S.
\endref

\bibitem{Ha} R. Hartshorne
\book Algebraic Geometry
\publ Springer-Verlag
\yr 1977
\endref

\bibitem{J} J-P. Jouanolou \book Theoremes de Bertini et applications.
\vol Prog.\ Math.\ 42 \publ Birkhauser \publaddr Boston \yr 1983
\endref

\bibitem{M} Y. I. Manin \book Cubic Forms \vol 4
\bookinfo North-Holland Mathematical Library
\yr 1986
\endref

\bibitem{Ma} A. Mayer \paper Families of K3 surfaces
\jour Nagoya J.\ Math.
\vol 48 \yr 1972
\pages 1--17
\endref

\bibitem{PS} I. Piateckii-Shapiro and I. R. Shafarevich
\paper A torelli theorem for
algebraic surfaces of type $K3$
\jour Math.\ USSR Izv. \vol 5 
\yr 1971
\pages 547--588
\endref

\bibitem{SD} B. Saint-Donat
\paper Projective models of $K3$ surfaces
\jour Amer.\ J.\ Math. \vol 96 
\yr 1974 \pages 602--639
\endref

\bibitem{Sk} F. Sakai
\paper Anticanonical models of rational surfaces
\jour Math.\ Ann. \vol 269 \yr 1984
\pages 389--410
\endref

\bibitem{St} H. Sterk
\paper Finiteness results for algebraic $K3$ surfaces
\jour Math.\ Z. \vol 189 \yr 1985 \pages 507--513
\endref

\bibitem{U} T. Urabe \paper On singularities on degenerate Del Pezzo surfaces of
degree $1$, $2$
\inbook Proc.\ Symp.\ Pure Math. \vol 40(2)
\yr 1983 \pages 587--592
\endref

\endRefs

\end{document}